\documentclass{amsart}
\usepackage[numbers, sort]{natbib}

\usepackage{amsmath}
\usepackage{amssymb}
\usepackage{amsmath}
\usepackage{amsfonts}
\usepackage{enumerate}
\usepackage{amsthm}
\usepackage{setspace}
\usepackage{multirow}
\usepackage{mathdots}
\usepackage{mathtools}
\usepackage{graphicx}
\usepackage{float}
\usepackage{mathrsfs}
\usepackage{cases}
\usepackage[all]{xy}
\usepackage{subcaption}

\usepackage{scalerel,stackengine}
\stackMath
\newcommand\reallywidehat[1]{%
\savestack{\tmpbox}{\stretchto{%
  \scaleto{%
    \scalerel*[\widthof{\ensuremath{#1}}]{\kern-.6pt\bigwedge\kern-.6pt}%
    {\rule[-\textheight/2]{1ex}{\textheight}}
  }{\textheight}%
}{0.5ex}}%
\stackon[1pt]{#1}{\tmpbox}%
}

\usepackage{natbib}

\DeclareMathOperator*{\ext}{ext}

\author{Xuefeng Shen}
\author{Melvin Leok}
\title[Lie group variational integrators using quaternions]{Lie group variational integrators for\\ rigid body problems using quaternions}

\begin{document}
\maketitle

\begin{abstract}
Rigid body dynamics on the rotation group have typically been represented in terms of rotation matrices, unit quaternions, or local coordinates, such as Euler angles. Due to the coordinate singularities associated with local coordinate charts, it is common in engineering applications to adopt the unit quaternion representation, and the numerical simulations typically impose the unit length condition using constraints or by normalization at each step. From the perspective of geometric structure-preserving, such approaches are undesirable as they are either computationally less efficient, or interfere with the preservation of other geometric properties of the dynamics.

In this paper, we adopt the approach used in constructing Lie group variational integrators for rigid body dynamics on the rotation group to the representation in terms of unit quaternions. In particular, the rigid body dynamics is lifted to unit quaternions, and the Lie group structure of unit quaternions is used to represent tangent vector intrinsically, thereby avoiding the use of a Lagrange multiplier. A Lie group variational integrator in the unit quaternion representation is derived, and numerical results are presented.
\end{abstract}

\section{Introduction}
For Lagrange mechanics on the tangent bundle $TQ$ and the Lagrangian $L$, the theory of variational integrators is well-established \cite{MaWe2001}, and is based on discretizing  Hamilton's principle rather than the differential equations of motion. The discrete Lagrangian $L_d: Q\times Q \to \mathbb{R}$ approximates the action integral over a small time interval,
\begin{align*}
    L_d(q_0, q_1) \approx \ext\limits_{q(0)=q_0, q(h)=q_1} \int_{0}^{h} L(q(t), \dot{q}(t)) dt,
\end{align*}
and the discrete Hamilton's principle states that the discrete action sum is stationary with respect to variations in the discrete solution curve that fixes the endpoints,
\begin{align*}
	\delta \sum_{k=0}^{N-1} L_d(q_k, q_{k+1})=0,	
\end{align*}
and this yields the discrete Euler Lagrange equations
\begin{align}\label{discreteEL}
    D_2Ld(q_0, q_1) + D_1(q_1, q_2) =0.
\end{align}
By introducing discrete Legendre transforms, we obtain a variational integrator on the cotangent bundle $T^*Q$,
\begin{equation}\label{typeI}
\begin{cases}
    p_0 =  -D_1L_d(q_0, q_1), \\
    p_1 =  D_2L_d(q_0, q_1).
\end{cases}
\end{equation}
It is easy to easy that \eqref{typeI} recovers \eqref{discreteEL} when the momentum variables are eliminated, and that \eqref{typeI} is the usual characterization of a symplectic map in terms of a Type I generating function. As such, variational integrators are automatically symplectic, and either form of the discrete equations can be easily implemented in a linear space or in local coordinates. However, the configuration manifold $Q$ is in general not a Euclidean space. Symplectic integrators do not exactly preserve energy, but instead rely on the existence of a modified energy associated for their long-time energy stability properties. Changing the coordinate chart at each time step results in the modified energy that is preserved changing from step to step, which destroys the long-time energy stability typically associated with symplectic integrators. In order to overcome this limitation, it is necessary to consider symplectic integrators that is are well-defined globally on the configuration manifold, and which commute with a change of coordinates.

One way to address this issue is to embed $Q$ in Euclidean space. In many cases, the configuration space naturally lives in a higher-dimensional Euclidean space $\mathbb{R}^n$, and is the level set of some constraint function $g: \mathbb{R}^n \to \mathbb{R}^m$. So the position and velocity are represented as vectors in the embedding space, and Lagrange multipliers are introduced to represent the constraint forces.

When configuration space $Q$ has a high codimension in the Euclidean space, and the constraint function is complicated, the above method is costly to implement. For example, the rotation group $SO(3)$, which is a three-dimensional Lie group, naturally lives in $\mathbb{R}^{3\times 3}$, a nine-dimensional Euclidean space, with six orthonormality constraints. As opposed to viewing $Q$ extrinsically as an embedded submanifold, when $Q$ is a Lie group, the group structure could be utilized to represent the position and tangent space intrinsically, thus avoiding the introduction of Lagrange multipliers. This idea was used in \cite{LeMcLe2005} for the rigid body problem, on configuration space $SO(3)$, and later generalized to the full body problem \cite{LeLeMc2007}.

Besides the rotation group $SO(3)$, the motion of rigid bodies can also be described by unit quaternions, which is an extremely popular approach in classical mechanics~\cite{Ho2011}, computer graphics~\cite{Go2009, Vi2011}, virtual reality~\cite{Ku1999}, and engineering applications in control~\cite{Ve2016, Ya2012} and estimation~\cite{Sh1978, Sh2006} of aerial and space vehicles. In \cite{WeMa1997}, a variational integrator for the rigid body problem was derived in terms of unit quaternions embedded in $\mathbb{R}^4$, so as to avoid the complicated constraints in $SO(3)$, but where the unit quaternion constraint is enforced using a Lagrange multiplier. We remark that the unit quaternions $S^3$, is not only the unit sphere in $\mathbb{R}^4$, but also a Lie group under quaternion multiplication.

In this paper, we will derive a Lie group variational integrator for rigid body problems using unit quaternions, but without the use of unit length constraints. Instead, we use the fact that $\mathbb{R}^3$ can be identified with the Lie algebra associated with $S^3\subset \mathbb{H}$. This will allow us to obtain a variational integrator for rigid body dynamics in the unit quaternion representation that is globally valid, expresses the relative rotation in terms of three degrees-of-freedom, and avoids the need for explicit constraints. In particular, it yields a numerical method in the unit quaternion representation that exhibits the computational advantages of Lie group variational integrators developed for rotation matrices, thereby making the advantages of geometric integrators accessible to engineering applications which are heavily invested in the unit quaternion representation of the rotation group.

\section{Background}
\subsection{Rigid body problem}
The configuration space of the rigid body is $SE(3)=\mathbb{R}^3\rtimes SO(3)$, where $(x, R)\in SE(3)$ describes the position and orientation of the body-fixed frame with respect to an inertial frame. Let $\rho\in \mathbb{R}^3$ denotes the position of mass element in the body-fixed frame, the corresponding inertial position is $x+R\rho$, and velocity $\dot{x}+\dot{R}\rho$, then the kinetic energy is given by
\begin{equation}\label{kinetic}
\begin{aligned}
T & = \frac{1}{2}\int_{\mathcal{B}}\|\dot{x}+\dot{R}\rho\|^2dm\\
&= \frac{1}{2}\int_{\mathcal{B}}\|\dot{x}\|^2dm + \int_{\mathcal{B}}\dot{x}^ \mathrm{T}\dot{R}\rho dm
        + \frac{1}{2}\int_{\mathcal{B}}\|\dot{R}\rho\|^2 dm\\
&= \frac{1}{2}m\|\dot{x}\|^2 + \frac{1}{2}\mathrm{tr}[\dot{R}J_d\dot{R}^ \mathrm{T}],
\end{aligned}
\end{equation}
where $J_d = \int_{\mathcal{B}}\rho\rho^{\mathrm{T}}dm$ is the nonstandard moment of inertia matrix. In \eqref{kinetic}, we used the fact that the origin of body-fixed frame is the center of mass of the body, thus
$$\int_{\mathcal{B}}\rho dm =0.$$
Since $R$ evolves on $SO(3)$, so $R^{\mathrm{T}}R = I$. By doing differentiation on both sides, we have $\dot{R}^{\mathrm{T}}R+R^{\mathrm{T}}\dot{R}=0$, thus $R^{\mathrm{T}}\dot{R}$ is skew symmetric. Define the hat map $\widehat{\cdot}: \mathbb{R}^3 \to R^{3\times 3}$ by the condition that $\widehat{x}y = x\times y$ for all $x,y\in \mathbb{R}^3$. If $x=(x_1,x_2,x_3)$, then $\widehat{x}$ is given by
\begin{displaymath}
  \widehat{x} = \begin{pmatrix}
0 &-x_3 &x_2\\
x_3 & 0 &-x_1\\
-x_2 & x_1 &0
\end{pmatrix},
\end{displaymath}
which is a skew-symmetric matrix. In particular, the hat map defines an isomorphism between $\mathbb{R}^3$ and skew-symmetric matrices. Since $R^{\mathrm{T}}\dot{R}$ is skew-symmetric, there exists a unique vector $\Omega$ such that $R^{\mathrm{T}}\dot{R}=\widehat{\Omega}$. In other words,
\begin{align}\label{rdot}
  \dot{R} = R\widehat{\Omega},
\end{align}
where $\Omega$ is the angular velocity in the body-fixed frame. This equation has a natural geometric interpretation: $SO(3)$ is a Lie group under matrix multiplication, and its Lie algebra, $\mathfrak{so}(3)$, consists of all skew-symmetric matrices,
\begin{align*}
  \mathfrak{so}(3) = \{\widehat{\Omega}\mid\Omega \in \mathbb{R}^3\}.
\end{align*}
The tangent space at each $R$ can be represented as the pushforward by the left action of $R$ on the Lie algebra,
\begin{align*}
  T_R SO(3) = (L_{R})_*\mathfrak{so}(3).
\end{align*}
In this way, the tangent bundle $TSO(3)$ is diffeomorphic to $SO(3)\times \mathfrak{so}(3)$ by left-trivialization. With this, we can express the rotational term in the kinetic energy as follows,
\begin{equation}\label{eq1}
  \begin{aligned}
    \mathrm{tr}[\dot{R}J_d\dot{R}^ \mathrm{T}] & = \mathrm{tr}[R \widehat{\Omega}J_d \widehat{\Omega}^{\mathrm{T}}R^\mathrm{T}] \\
     & = \mathrm{tr}[\widehat{\Omega}J_d \widehat{\Omega}^{\mathrm{T}}]\\
     & = \Omega^\mathrm{T}J\Omega,
  \end{aligned}
\end{equation}
where $J = \mathrm{tr}[J_d]I_3 - J_d$ is standard moment of inertia matrix. Combining \eqref{kinetic} and \eqref{eq1}, we obtain $T = \frac{1}{2}m\|\dot{x}\|^2+ \frac{1}{2}\Omega^\mathrm{T}J\Omega$, and the Lagrangian for rigid body problems $L: TSE(3) \to R$ is given by,
\begin{align*}
  L(x,R,\dot{x},\Omega) = T-V = \frac{1}{2}m\|\dot{x}\|^2+ \frac{1}{2}\Omega^\mathrm{T}J\Omega - V(x,R).
\end{align*}

\subsection{Unit quaternions}
Quaternions are a noncommutative division algebra, usually denoted as $\mathbb{H}$, with element $q=(q_s, \vec{q}_v)$, where $q_s$ is the scalar part and $\vec{q}_v \in \mathbb{R}^3$ is the vector part.
$\mathbb{H}$ inherits the vector space structure and differential structure from $\mathbb{R}^4$, and has the following multiplication multiplication: given $q=(q_s, \vec{q}_v)$, $p=(p_s, \vec{p}_v)$,
\begin{align*}
  q\cdot p = (q_s\cdot p_s - \vec{q}_v\cdot \vec{p}_v, q_s\cdot \vec{p}_v+ p_s\cdot \vec{q}_v+ \vec{q}_v\times \vec{p}_v).
\end{align*}
We can also define conjugatation $q^* = (q_s, -\vec{q}_v)$, norm $\| q\| = \sqrt{qq^*}$, and inverse under quaternion multiplication $q^{-1} = \frac{q^*}{\| q\|^2}$. The set of unit quaternions
$$S^3 = \{q\in \mathbb{H}\mid \| q\|=1\}$$
is diffeomorphic to the unit sphere in $\mathbb{R}^4$, and forms a Lie group under quaternion multiplication. The identity element is $e = (1,0,0,0)$, and the Lie algebra $$T_e S^3 = \{(0, \xi)\mid \xi\in \mathbb{R}^3\}.$$
The corresponding exponential map is
\begin{align}
\exp(\xi) = \bigg(\mathrm{cos}(|\xi|),\frac{\xi}{|\xi|}\mathrm{sin}(|\xi|)\bigg).
\end{align}
For any $q_0\in S^3$, the map $\xi\mapsto q_0\cdot \exp(\xi)$ gives a local diffeomorphism between a neighborhood of the origin of $\mathbb{R}^3$ and a neighborhood of $q_0$. It can be verified that for $q\in S^3, v\in \mathbb{R}^3$, $q(0,v)q^*$ is pure imaginary, i.e., the scalar part vanishes. This defines a linear action on the vector part, which turns out to be a rotation. Thus, for each $q\in S^3$, we have the corresponding rotation $\pi(q)\in SO(3)$,
\begin{align*}
  S^3 \xrightarrow{\mathmakebox[1cm]{\pi}} SO(3),\qquad \pi(q) = (2q^2_s -1)I_3 +2\vec{q}_v\vec{q}^\mathrm{T}_v +2q_s \widehat{\vec{q}_v}.
\end{align*}
The map $\pi$ is surjective, locally diffeomorphic, and is also a Lie group homomorphism:
$$\pi(q_1\cdot q_2) = \pi(q_1)\cdot \pi(q_2),$$ where the operation on left side is quaternion multiplication, and the operation on the right side is matrix multiplication. $\pi$ is not a global diffeomorphism, for each $R\in SO(3)$, its preimage is always one pair of antipodal points $\pm q$. Actually, $S^3$ and $SO(3)$ are not diffeomorphic, as they have different homology groups.

We can lift the Lagrangian dynamics on $TSE(3)$ to $T\mathbb{R}^3\times TS^3$:
$$T\mathbb{R}^3\times TS^3 \xrightarrow{\mathmakebox[1cm]{Id\times T\pi}} TSE(3)\xrightarrow{\mathmakebox[1cm]{L}} R,$$
As manifolds, we have $TSE(3) = T(\mathbb{R}^3\times SO(3))=T\mathbb{R}^3\times TSO(3)$, where $Id$ is the identity map from $T\mathbb{R}^3$ to itself, and $T\pi: TS^3\to TSO(3)$ is the tangent lift of $\pi$. Define the Lagrangian $\hat{L}= L\circ (Id\times T\pi)$ on $T(\mathbb{R}^3\times S^3)$, and we will develop our algorithm for the Lagrangian mechanics defined by $\hat{L}$ on $T(\mathbb{R}^3\times S^3)$. Recall that $Id\times T\pi$ is a local diffeomorphism, so $(T(\mathbb{R}^3\times S^3), \hat{L})$ and $(TSE(3), L)$ are equivalent for initial-value problems, and locally equivalent for two point boundary-value problems. A calculation shows that,
$$(q, q\cdot (0, \xi))\xrightarrow{\mathmakebox[1cm]{T\pi}} (\pi(q), \pi(q)\cdot \reallywidehat{2\xi}),$$
so
\begin{equation}\label{lhat}
\begin{aligned}
  \hat{L}(x,q,\dot{x},\dot{q}) & =  \hat{L}(x,q,\dot{x},q\cdot(0,\xi))\\
                               & =  \hat{L}(x,q,\dot{x},\xi)\\
                               & = \frac{1}{2}m\|\dot{x}\|^2+ 2\xi^\mathrm{T}J\xi - V(x,q).
\end{aligned}
\end{equation}

\section{Lie group variational integrator}
\subsection{Continuous time equation}
We derive Euler--Lagrange equations for the Lagrangian $\hat{L}$ on $T(\mathbb{R}^3\times S^3)$. By Hamilton's principle, $(x(t),q(t))\in \mathbb{R}^3\times S^3$ extremizes the action integral
$$\int_{t_0}^{t_1}\hat{L}(x(t),q(t),\dot{x}(t),\dot{q}(t))dt,$$
for variations that fix the endpoints $x(t_0)=x_0, x(t_1)=x_1$ and $q(t_0)=q_0, q(t_1)=q_1$.
Consider variations of $(x(t),q(t))$ parameterized by $\lambda$: Given any $(\delta x(t), \eta(t))$ that vanish at the endpoints, we construct $(x(t,\lambda), q(t, \lambda))$ as follows,
\begin{align*}
  x(t,\lambda) = x(t) + \lambda \delta x(t),
\end{align*}
and $q(t,\lambda)\in S^3$, such that $q(t,0)=q(t)$, satisfies the following equation,
\begin{align}\label{deltaepsilon}
	\frac{\partial q}{\partial \lambda}(t,0) &= q(t)\cdot (0, \eta(t)),
\end{align}
which implies that $q(t_0,\lambda)= q_0$ and $q(t_1,\lambda)= q_1$, since $\eta(t)$ vanishes at the endpoints. Since $\frac{\partial q}{\partial t}(t, \lambda)\in T_{q(t,\lambda)}S^3$, we have by left-trivialization, $\frac{\partial q}{\partial t}(t, \lambda)= q(t,\lambda)\cdot (0, \xi(t,\lambda))$, for a suitable choice of $\xi(t,\lambda)\in\mathbb{R}^3$. Taking derivatives with respect to $\lambda$ on both sides, we obtain
\begin{align*}
  \frac{\partial^2 q}{\partial \lambda\partial t}(t,\lambda) = \frac{\partial q}{\partial \lambda}(t, \lambda)\cdot(0, \xi(t,\lambda)) +q(t,\lambda)\cdot \bigg(0, \frac{\partial \xi}{\partial \lambda}(t,\lambda)\bigg).
\end{align*}
Evaluating this at $\lambda=0$ yields,
\begin{equation}\label{eq3}
\begin{aligned}
  \frac{\partial^2 q}{\partial \lambda\partial t}(t,0) &= \frac{\partial q}{\partial \lambda}(t, 0)\cdot(0, \xi(t)) +q(t)\cdot \bigg(0, \frac{\partial \xi}{\partial \lambda}(t,0)\bigg)\\
&= q(t)\cdot (0, \eta(t))\cdot(0, \xi(t)) +q(t)\cdot \bigg(0, \frac{\partial \xi}{\partial \lambda}(t,0)\bigg),
\end{aligned}
\end{equation}
where we used \eqref{deltaepsilon}, and we let $\xi(t)=\xi(t,0)$. Taking derivatives with respect to $t$ on both sides of \eqref{deltaepsilon} yields,
\begin{equation}\label{eq4}
  \begin{aligned}
    \frac{\partial^2 q}{\partial t\partial \lambda}(t,0) & = \dot{q}(t)\cdot (0, \eta(t))+ q(t)\cdot(0,\dot{\eta}(t)) \\
    & = q(t)\cdot (0,\xi(t))\cdot (0, \eta(t))+ q(t)\cdot(0,\dot{\eta}(t)).
  \end{aligned}
\end{equation}
Equating \eqref{eq3} and \eqref{eq4} by the equality of mixed partials, we get
\begin{equation*}
\begin{aligned}
  \bigg(0,\frac{\partial \xi}{\partial \lambda}(t,0)\bigg) &=  (0,\dot{\eta}(t))+ (0,\xi(t))\cdot (0, \eta(t))-(0, \eta(t))\cdot(0, \xi(t))\\
  &= (0,\dot{\eta}(t))+ (0, 2\xi(t)\times \eta(t)).
  \end{aligned}
\end{equation*}
Thus,
\begin{align}\label{eq5}
  \frac{\partial \xi}{\partial \lambda}(t,0) = \dot{\eta}(t) + 2\xi(t)\times \eta(t).
\end{align}

Using the Lagrangian given in \eqref{lhat}, Hamilton's principle states that
\begin{align*}
  \left.\frac{d}{d\lambda}\right|_{\lambda=0}\int_{t_0}^{t_1}\bigg[\frac{1}{2}m\| \dot{x}+\lambda \dot{\delta x}\|^2 + 2\xi(t,\lambda)^\mathrm{T}J\xi(t,\lambda)-V(x(t,\lambda),q(t,\lambda))\bigg]dt =0,
\end{align*}
which means
\begin{align*}
  \int_{t_0}^{t_1} \bigg[m\dot{x}\cdot \dot{\delta x} +4\xi(t)^\mathrm{T}J\frac{\partial \xi}{\partial \lambda}(t,0)-\bigg(\frac{\partial V}{\partial x}\cdot \delta x + \bigg(\frac{\partial V}{\partial q}\bigg)^\mathrm{T}\frac{\partial q}{\partial \lambda}(t,0)\bigg)\bigg] dt =0.
\end{align*}
The terms involving the infinitesimal variation $\delta x$ are
\begin{align}\label{eqdeltax}
  \int_{t_0}^{t_1}\bigg[m\dot{x}\cdot \dot{\delta x} - \frac{\partial V}{\partial x}\cdot \delta x\bigg] dt = \int_{t_0}^{t_1}\bigg(-m\ddot{x}-\frac{\partial V}{\partial x}\bigg)\delta x dt,
\end{align}
and using \eqref{eq5}, we have that
\begin{equation}
\begin{aligned}
\int_{t_0}^{t_1} 4\xi(t)^\mathrm{T}J\frac{\partial \xi}{\partial \lambda}(t,0) dt &= \int_{t_0}^{t_1} 4\xi^\mathrm{T}J(\dot{\eta} + 2\xi\times \eta) dt\\
& = \int_{t_0}^{t_1} 4(-\dot{\xi}^\mathrm{T}J + 2 \xi^\mathrm{T}J\widehat{\xi})\eta dt.
\end{aligned}
\end{equation}

For any $q\in \mathbb{H}, v\in \mathbb{R}^3$, define $F(q): \mathbb{H}\to \mathbb{R}^{3\times 4}$ by the condition that
$$q\cdot (0,v) = F(q)^\mathrm{T}v.$$
It can be easily verified that for $q=(q_s, \vec{q}_v)$, $F(q)=(-\vec{q}_v, q_s\cdot I - \reallywidehat{\vec{q}_v})$, thus
\begin{equation}\label{eq6}
\begin{aligned}
\int_{t_0}^{t_1}\bigg(\frac{\partial V}{\partial q}\bigg)^\mathrm{T}\frac{\partial q}{\partial \lambda}(t,0)dt &= \int_{t_0}^{t_1}\bigg(\frac{\partial V}{\partial q}\bigg)^\mathrm{T}(q\cdot(0,\eta)) dt\\
& = \int_{t_0}^{t_1}\bigg(\frac{\partial V}{\partial q}\bigg)^\mathrm{T}F(q)^\mathrm{T}\eta dt.
\end{aligned}
\end{equation}
Combining \eqref{eq3}, \eqref{eqdeltax}, \eqref{eq6}, and integrating by parts, we have
\begin{align*}
 \int_{t_0}^{t_1} \bigg(-m\ddot{x}-\frac{\partial V}{\partial x}\bigg)\delta x +\bigg(-4\dot{\xi}^\mathrm{T}J+8\xi^\mathrm{T}J\widehat{\xi}+\bigg(\frac{\partial V}{\partial q}\bigg)^\mathrm{T}F(q)^\mathrm{T}\bigg)\eta dt =0,
\end{align*}
for all variations $\delta x$ and $\eta$ that vanish at the endpoints. By the fundamental theorem of the calculus of variations, the Euler--Lagrange equations for the Lagrangian $\hat{L}$ on $T(\mathbb{R}^3\times S^3)$ is given by,
\begin{equation}
\left\{
\begin{aligned}
m\ddot{x} &= -\frac{\partial V}{\partial x}, \\
4J\dot{\xi} + 8\xi\times (J\epsilon) &= F(q)\frac{\partial V}{\partial q},\\
\dot{q} &= q\cdot (0,\xi).
\end{aligned}\right.
\end{equation}

\subsection{Variational integrator on the Lagrangian side}
The discrete Lagrangian from \cite{WeMa1997} is used here, which can be viewed as a midpoint rule approximation of the integral, combined with linear interpolation. Given endpoints $(x_0, q_0)$ and $(x_1, q_1)$, since $(0,\xi)=q^* \dot{q}\approx (\frac{q_0+q_1}{2})^*\cdot(\frac{q_1-q_0}{h})= (0, \frac{1}{h}\text{Im}(q^*_0q_1))$, we have that $\xi$ is approximated by $\frac{1}{h}\text{Im}(q^*_0q_1)$. We can construct
the discrete Lagrangian as
\begin{equation}\label{discrtelag}
\begin{aligned}
  &L_d(x_0,q_0,x_1,q_1)\\
  &\qquad=h\bigg( \frac{1}{2}m\| \frac{x_1-x_0}{h}\|^2 + \frac{2}{h^2}(\text{Im}(q^*_0q_1))^\mathrm{T}J(\text{Im}(q^*_0q_1))-\frac{V(x_0,q_0)+V(x_1,q_1)}{2}\bigg).
\end{aligned}
\end{equation}
The discrete Euler Lagrange equation is
\begin{numcases}\\
D_{x_1}(L_d(x_0,q_0,x_1,q_1)+L_d(x_1,q_1,x_2,q_2)) =0,  \label{del1}\\
D_{q_1}(L_d(x_0,q_0,x_1,q_1)+L_d(x_1,q_1,x_2,q_2)) =0.   \label{del2}
\end{numcases}
By substituting the expression for the discrete Lagrangian \eqref{discrtelag} into \eqref{del1}, we obtain
\begin{align}\label{del11}
  m\cdot \frac{x_2-2x_1+x_0}{h^2} = -\frac{\partial V}{\partial x}(x_1,q_1).
\end{align}
For \eqref{del2}, recall that $q_1$ evolves on $S^3$, so we consider a variation $q_1(\lambda)$ of $q_1$, such that $q_1(0)=q_1$ and $\delta{q}_1=\left.\frac{d q_1(\lambda)}{d \lambda}\right|_{\lambda=0}=q_1\cdot(0,\eta)$, then
\begin{equation}\label{derivq0q1}
	\begin{aligned}
  	\left.\frac{d}{dt}\right|_{\lambda=0}(\text{Im}(q^*_0q_1(\lambda)))^\mathrm{T}J(\text{Im}(q^*_0q_1(\lambda)))
  	&= 2(\text{Im}(q^*_0q_1))^\mathrm{T}J\bigg(\!\left.\frac{d}{d\lambda}\right|_{\lambda=0}\text{Im}(q^*_0q_1(\lambda))\bigg)\\
   	&= 2(\text{Im}(q^*_0q_1))^\mathrm{T}J(\text{Im}(q^*_0 q_1(0,\eta))).
   \end{aligned}
\end{equation}
Define $G(q): \mathbb{H}\to \mathbb{R}^{3\times 3}$, such that for any $q\in \mathbb{H}$, $v\in \mathbb{R}^3$,
\begin{align*}
  \text{Im}(q\cdot(0,v)) = G(q)^\mathrm{T}v.
\end{align*}
It can be verified that
$$G(q) = q_s\cdot I_3 - \reallywidehat{\vec{q}_v},$$
so $$2(\text{Im}(q^*_0q_1))^\mathrm{T}J(\text{Im}(q^*_0 q_1(0,\eta))) = 2(\text{Im}(q^*_0q_1))^\mathrm{T}JG(q^*_0q_1)^\mathrm{T}\eta.$$
Similarly,
\begin{equation}\label{derivq1q2}
	\begin{aligned}
  		\left.\frac{d}{dt}\right|_{\lambda=0}(\text{Im}(q^*_1(\lambda)q_2))^\mathrm{T}J(\text{Im}(q^*_1(\lambda)q_2)) & = 2(\text{Im}(q^*_1q_2))^\mathrm{T}J(\text{Im}((q_1(0,\eta))^* q_2)\\
    	& = 2(\text{Im}(q^*_2q_1))^\mathrm{T}J(\text{Im}(q^*_2 q_1(0,\eta)))\\
    	& = 2(\text{Im}(q^*_2q_1))^\mathrm{T}JG(q^*_2 q_1)^\mathrm{T}\eta.
    \end{aligned}
\end{equation}
Substituting \eqref{derivq0q1} and \eqref{derivq1q2} into \eqref{del2} yields
$$\frac{4}{h^2}(\text{Im}(q^*_0q_1))^\mathrm{T}JG(q^*_0q_1)^\mathrm{T}\eta + \frac{4}{h^2}(\text{Im}(q^*_2q_1))^\mathrm{T}JG(q^*_2q_1)^\mathrm{T}\eta -\bigg(\frac{\partial V}{\partial q}(x_1, q_1)\bigg)^\mathrm{T}F(q_1)^\mathrm{T}\eta = 0,$$
which holds for any $\eta$, which gives
\begin{align}\label{del22}
  G(q^*_0q_1)J(\text{Im}(q^*_0q_1)) + G(q^*_2q_1)J(\text{Im}(q^*_2q_1)) = \frac{h^2}{4}F(q_1)\frac{\partial V}{\partial q}(x_1, q_1).
\end{align}
In summary, the discrete Euler--Lagrange equations derived from the discrete Lagrangian \eqref{discrtelag} are given by
\begin{equation*}
\left\{
  \begin{aligned}
m\cdot \frac{x_2-2x_1+x_0}{h^2} &= -\frac{\partial V}{\partial x}(x_1,q_1), \\
G(q^*_0q_1)J(\text{Im}(q^*_0q_1)) + G(q^*_2q_1)J(\text{Im}(q^*_2q_1)) &= \frac{h^2}{4}F(q_1)\frac{\partial V}{\partial q}(x_1, q_1),
  \end{aligned}\right.
\end{equation*}
where $G(q) = q_s\cdot I_3 - \reallywidehat{\vec{q}_v}$, $F(q)=(-\vec{q}_v, q_s\cdot I - \reallywidehat{\vec{q}_v})$.

\subsection{Variational integrator on the Hamiltonian side}
For any Lie group $G$, its tangent bundle $TG$ is diffeomorphic to $G\times \mathfrak{g}$ by left-trivialization,
\begin{alignat*}{2}
	TG &\cong G\times \mathfrak{g},  &\quad (q,v_q)&\xrightarrow{\mathmakebox[1cm]{\psi}}(q,L_{q^{-1}*}v_q),
\intertext{and $T^* G$ is diffeomorphic to $G\times \mathfrak{g}^*$,}
T^* G &\cong G\times \mathfrak{g}^*, &\quad  (q, p_q)&\xrightarrow{\mathmakebox[1cm]{\phi}} (q, {L_q}^*p_q).
\end{alignat*}
Given a Lagrangian $L: TG\to R$, we can define $\hat{L}=L\circ \psi^{-1}$, and the following diagram commutes,
\begin{equation*}
    \xymatrix@!0@R=0.75in@C=1in{
      TG \ar[r]^{F_L} \ar[d]_{\psi}& T^{*}G \ar[d]^{\phi}
      \\
      G\times \mathfrak{g}  \ar[r]_{F_{\hat{L}}} & G\times \mathfrak{g}^*
    }
\end{equation*}
where $F_L$ is the Legendre transformation and $F_{\hat{L}}$ denotes the partial derivative of $\hat{L}$ with respect to the Lie algebra.

For a Lagrangian of the form $\hat{L}(x,q,\dot{x},\xi) = \frac{1}{2}m\|\dot{x}\|^2+ 2\xi^\mathrm{T}J\xi - V(x,q)$, the Legendre transform is given by
$(x,q,\dot{x},\xi)\mapsto (x,q, m\dot{x}, 4J\xi)$.

For a discrete Lagrangian $L_d: G\times G\to R$, we consider a variation $q_0(\lambda)$ of $q_0$, such that $q_0(0) = q_0$ and $\delta{q_0}=\left.\frac{d q_0(\lambda)}{d \lambda}\right|_{\lambda=0}= L_{q_0*}\xi$, then an element $\alpha\in \mathfrak{g}^*$ is defined by the condition
$$\left.\frac{d}{d\lambda}\right|_{\lambda=0}L_d(q_0(\lambda),q_1) = \langle\alpha, \xi\rangle.$$
It is easy to verify that $\alpha = {L_{q_0}}^*D_1L_d(q_0, q_1)$. The same arguments hold for the derivative with respect to $q_1$. This allows us to conclude that the the partial derivative of the discrete Lagrangian with respect to the Lie algebra element in the left-trivialized coordinate system is related to the partial derivative with respect to the Lie algebra element that generates the variation of $q_i\in G$ is related to the partial derivative with respect to $q_i$ by a left translation. So the variational integrator on the Hamiltonian side for the discrete Lagrangian \eqref{discrtelag} is given by
\begin{equation}\label{eqhamilton}
 \left\{
\begin{aligned}
 p_0 &= m\cdot \frac{x_1-x_0}{h} +\frac{h}{2}\frac{\partial V}{\partial x}(x_0, q_0),  \\
 w_0 &= -\frac{4}{h}G(q^*_1q_0)J(Im(q^*_1q_0)) + \frac{h}{2}F(q_0)\frac{\partial V}{\partial q}(x_0,q_0), \\
 p_1 &= m\cdot \frac{x_1-x_0}{h} -\frac{h}{2}\frac{\partial V}{\partial x}(x_1, q_1), \\
 w_1 &= \frac{4}{h}G(q^*_0q_1)J(Im(q^*_0q_1)) - \frac{h}{2}F(q_1)\frac{\partial V}{\partial q}(x_1,q_1).
\end{aligned}\right.
\end{equation}
Here, $w_0, w_1$ are in the dual space of $\mathfrak{so}(3)$, and the derivation follows the approach adopted in the previous section.

\section{Implementation of algorithm}
Given initial conditions $(x_0, R_0, \dot{x}_0, \dot{R}_0=R_0 \widehat{\Omega_0})$ for the rigid body problem, our algorithm is given as follows,
\begin{enumerate}
\item lift $(x_0, R_0)\in SE(3)$ to $(x_0, q_0)\in \mathbb{R}^3\times S^3$, such that $\pi(q_0) = R_0$ (see~\cite{La2008});
\item start with initial conditions $(x_0, q_0, p_0 = m\dot{x}_0, w_0 = 2J\Omega_0)$ on $T^*(\mathbb{R}^3\times S^3)$;
\item apply \eqref{eqhamilton} repeatedly to generate the discrete sequence $(x_k, q_k, p_k, w_k)$.

\end{enumerate}
At each iteration, $x_1$ can be directly calculated by
$$x_1 = x_0 +\frac{h}{m}(p_0-\frac{h}{2}\frac{\partial V}{\partial x}(x_0,q_0)).$$
For $q_1$, we represent it via exponential map $$q_1=q_0\cdot \exp(\xi_0),$$ for an element $\xi_0\in\mathbb{R}^3$. The root finding problem of determining the $\xi_0$ that satisfies the second equation in the nonlinear system \eqref{eqhamilton} can be performed by approximating the Jacobian matrix by a finite-difference approximation at the beginning, and then applying Broyden's method to update the Jacobian matrix at each step, which avoids the complicated task for evaluating the Jacobi matrix exactly. Then, $p_1, w_1$ can both be calculated explicitly from the last two equations in the system \eqref{eqhamilton}. If we wish to recover the discrete solution in the rotation matrix representation, we can  project $(x_k, q_k, p_k, w_k)$ to $TSE(3)$ by using
\begin{align*}
R_k &= \pi(q_k),\\
\dot{x}_k &= \frac{p_k}{m},\\
\Omega_k &= \frac{1}{2}J^{-1}w_k.	
\end{align*}

\section{Numerical Experiment}
\begin{figure}
\includegraphics[scale=0.4]{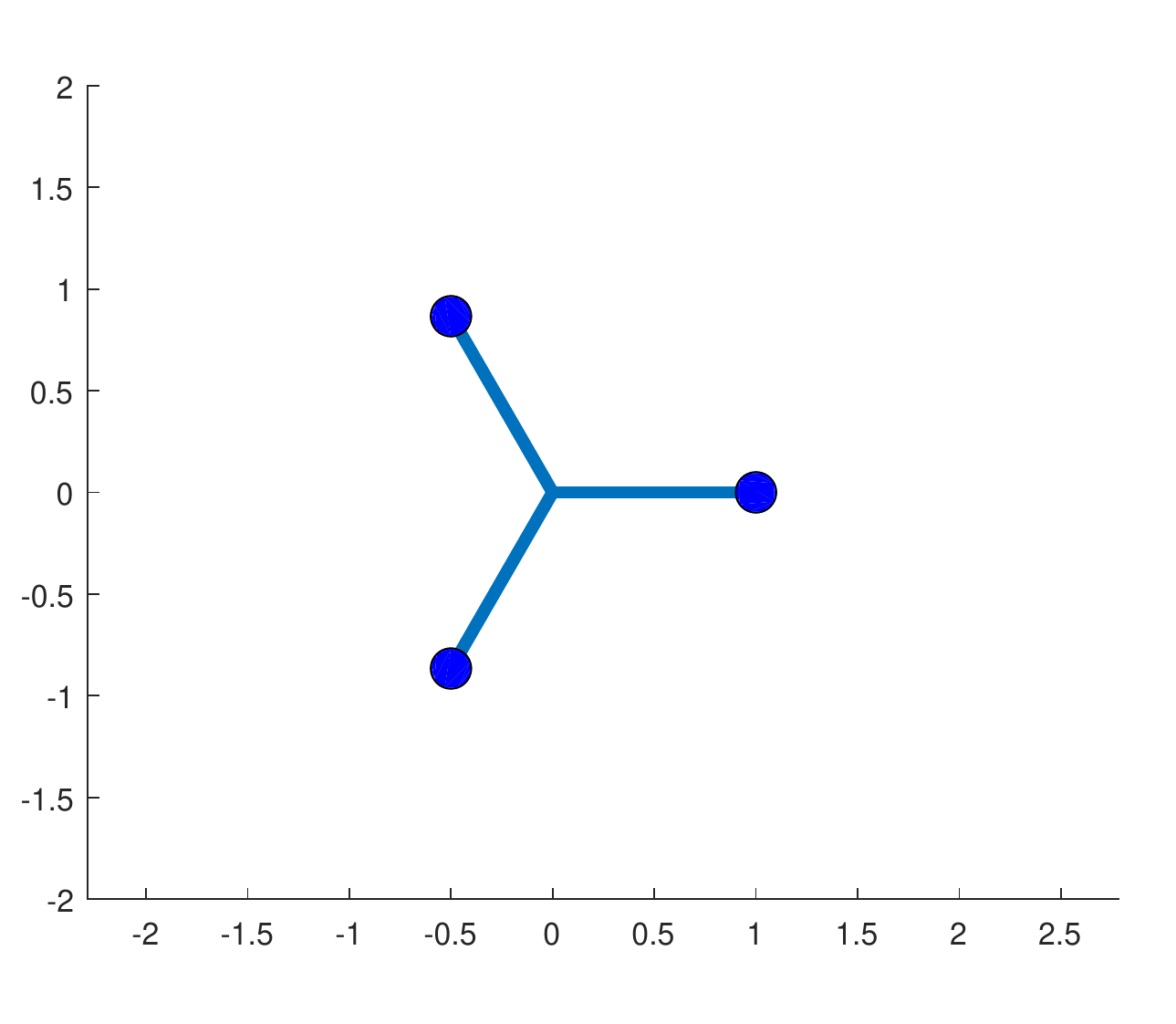}
\caption{Planar rigid body\label{fig:planar_body}}
\end{figure}

\begin{figure}
	\begin{subfigure}[b]{0.45\textwidth}
		\includegraphics[scale=0.4]{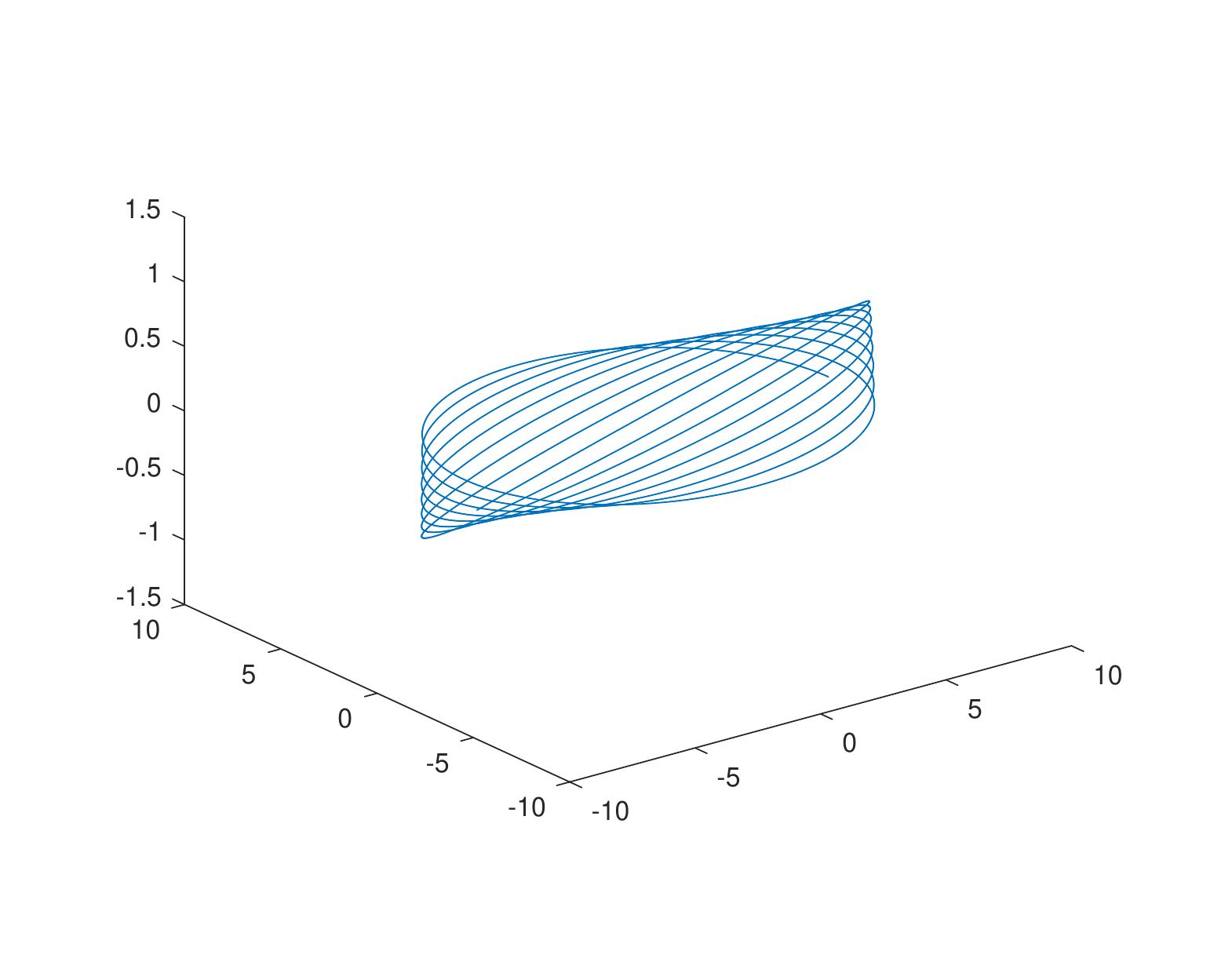}
		\caption{trajectory of the center of mass}
	\end{subfigure}
	\begin{subfigure}[b]{0.45\textwidth}
		\includegraphics[scale=0.4]{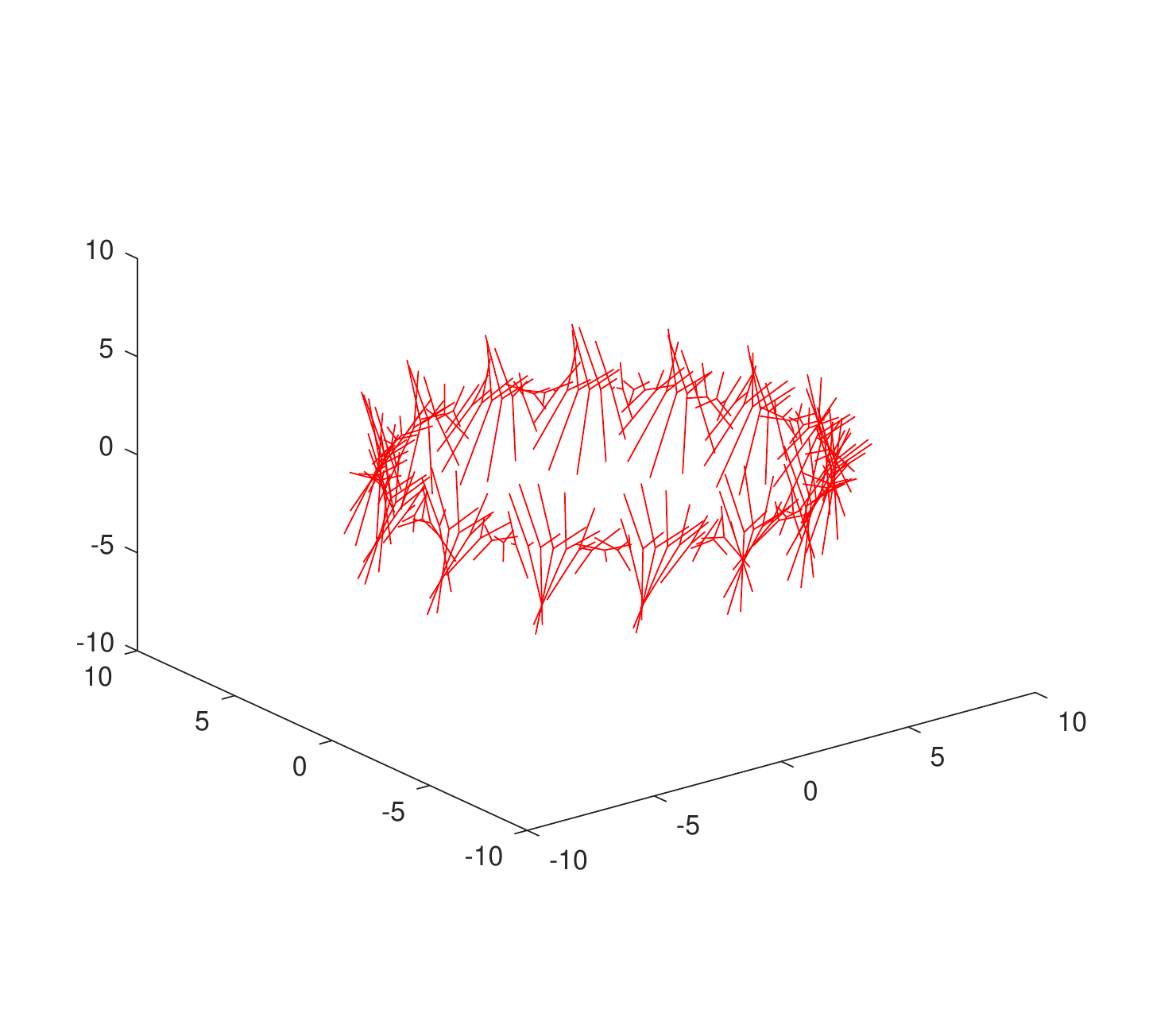}
		\caption{rigid body motion over one orbit}
\end{subfigure}
\caption{Motion of the rigid body in a central gravitational field.\label{fig:trajectories}}
\end{figure}

\begin{figure}
	\begin{subfigure}[b]{0.45\textwidth}
		\includegraphics[scale=0.4]{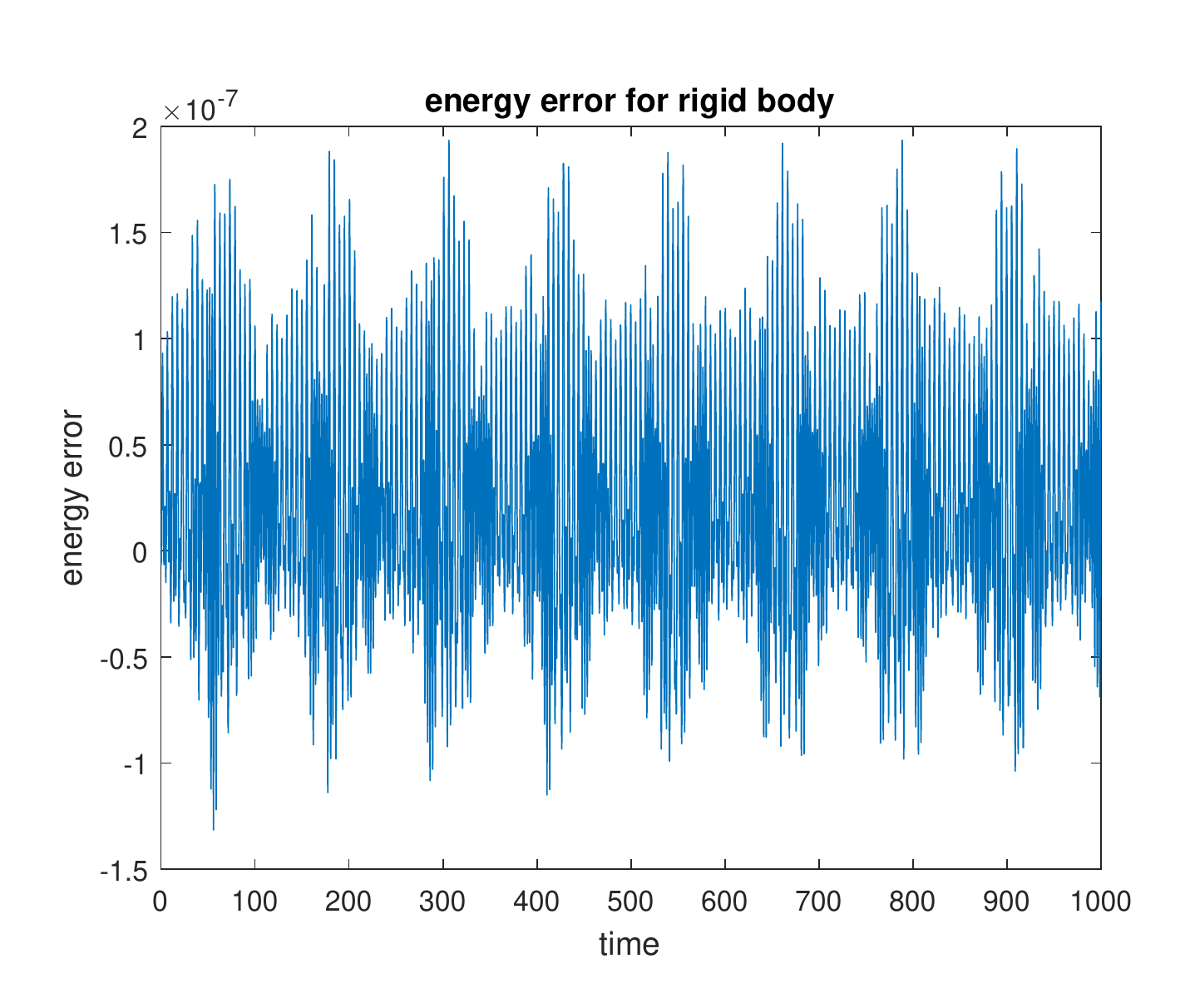}
		\caption{energy error}
	\end{subfigure}
	\begin{subfigure}[b]{0.45\textwidth}
		\includegraphics[scale=0.4]{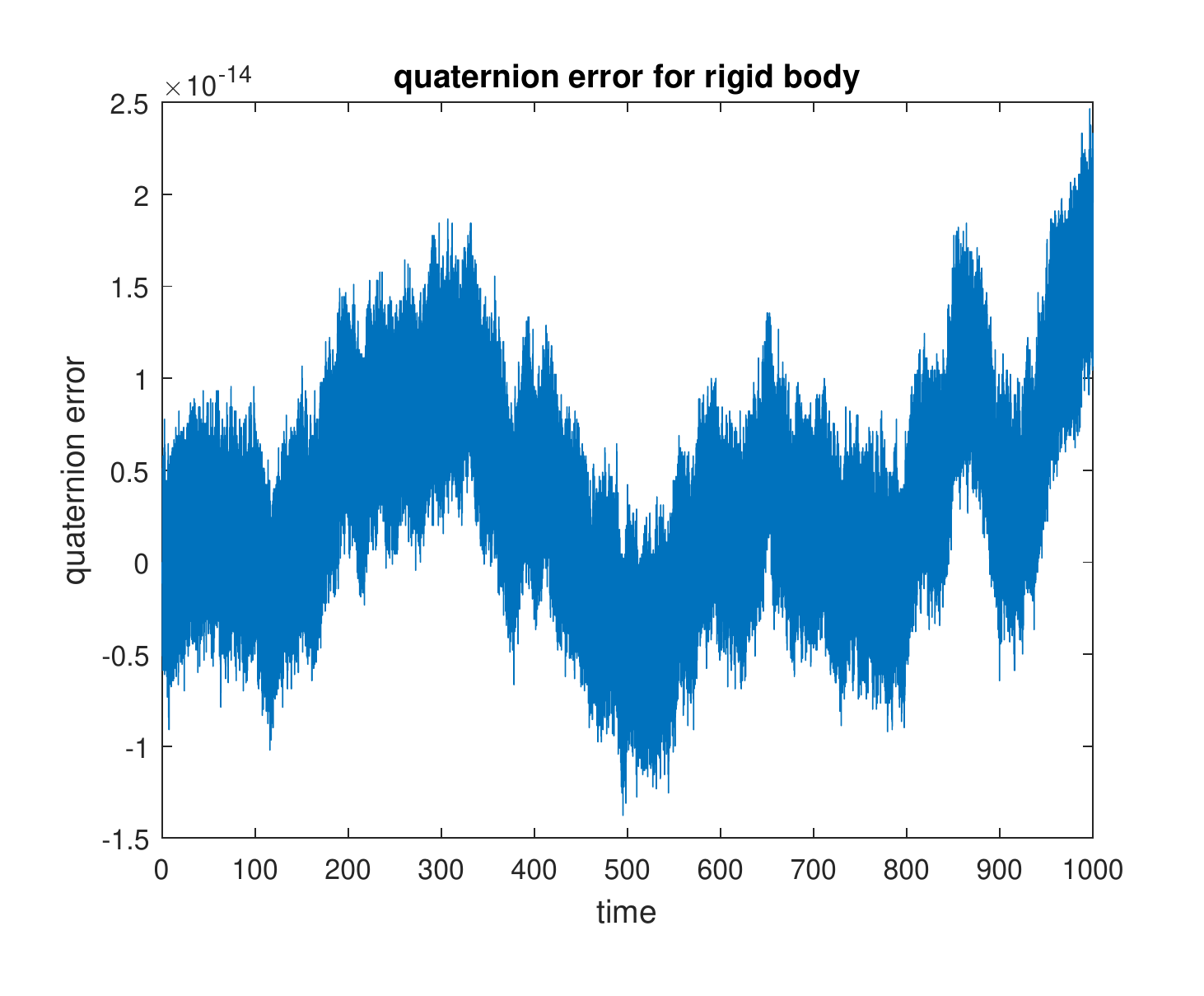}
		\caption{quaternion error}
\end{subfigure}
\caption{The Lie group variational integrator using quaternions preserves the energy and the unit norm condition very well.\label{fig:errors}}
\end{figure}

\begin{figure}
\includegraphics[scale=0.4]{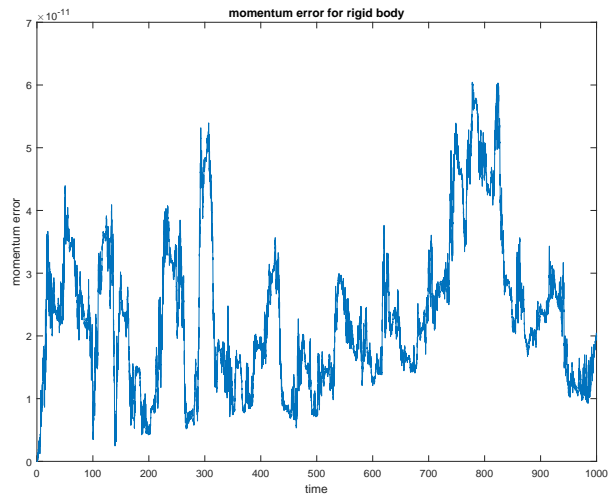}
\caption{Total angular momentum error. \label{fig:momentum}}
\end{figure}

We consider a planar rigid body composed of three uniform balls with unit mass and radius 0.1, connected by massless rods as in Figure~\ref{fig:planar_body}.
We consider the motion of the planar rigid body in a gravitational field centered at the origin with potential energy $-\frac{1}{r}$, with initial conditions $x_0 =(8,0,0) $, $q_0 = (1,0,0,0)$, $p_0 = (1,0,0)$, and $w_0 = (1,2,3)$, and timestep $h=0.01$. The trajectory of the center of mass during the time interval $[0,1000]$ and motion of the rigid body over one orbital period are given in Figure~\ref{fig:trajectories}.
The energy error and quaternion error, which is the deviation from the unit quaternion surface, during time $[0,1000]$ are given in Figure~\ref{fig:errors}, and we observe that these invariants are preserved extremely well.
It can be seen that energy error remains stable due to the symplectic property of algorithm, and quaternion error remains stable due to the intrinsic expression of Lie group. Also, we observe that the Lagrangian \eqref{lhat} is invariant under the group action $S^3\times (R^3\times S^3)\mapsto R^3\times S^3$ given by
$$q_1(x,q) = (\pi(q_1)x, q_1q).$$
Thus, the corresponding momentum map, which is the total angular momentum $x\times p + \pi(q)J\Omega$, is preserved. The discrete Lagrangian \eqref{discrtelag} we designed is also invariant under the diagonal group action, so by the discrete Noether's theorem \cite{MaWe2001}, our algorithm also preserves total angular momentum of the system. The total angular momentum 
error is given in Figure \ref{fig:momentum}. As observed, it is preserved quite well.

\section{Conclusions}
A Lie group type variational integrator is developed for rigid body dynamics represented on the unit quaternion group. The continuous equations of motion, and variational integrators on both the Lagrangian and Hamiltonian sides are derived. The resulting algorithm is symplectic, and exhibits good long time stability of the energy and the unit quaternion condition.

By representing the dynamics in terms of unit quaternions as opposed to the rotation group, one requires less memory to store the attitude, which can be valuable when considering the dynamics of a large collection of rigid bodies. Furthermore, since quaternion multiplication is cheaper than the multiplication of rotation matrices, this could also potentially result in reduced computational costs. In particular, the memory requirements for quaternions is $4$ scalars, and for $3\times 3$ matrices is $9$ scalars, and the cost of composing quaternions is $16$ multiplies, and the cost of composing $3\times 3$ matrices is $27$ multiplies.

Most importantly, many applications in computer graphics, and guidance, navigation, control, and estimation of aerial and space vehicles relies on the representation of the attitude in terms of unit quaternions, so the construction of Lie group variational integrators using the unit quaternion representation allows users with an investment in terms of representing rotations using unit quaternions can now take advantage of the benefits of Lie group variational integrator techniques, such as symplecticity, long-time energy stability, and the ability to respect the unit quaternion condition without the use of local coordinates or constraints.

\section*{Acknowledgements} The authors were supported in part by NSF grants DMS-1001521, CMMI-1029445, DMS-1065972, CMMI-1334759, DMS-1411792, and NSF CAREER award DMS-1010687.

\bibliography{quaternion}
\bibliographystyle{plainnat}

\end{document}